\documentclass[authoryear]{elsarticle}

\usepackage{graphicx}
\usepackage{epstopdf}
\usepackage{amsfonts}
\usepackage{color}
\usepackage{mathrsfs}
\usepackage{psfrag, amssymb, amsmath, cite}
\usepackage{verbatim}
\usepackage{epsfig} 
\usepackage{times} 
\usepackage{amsmath} 
\usepackage{amssymb}  
\usepackage{amsfonts}
\usepackage[colorlinks,linkcolor=red,anchorcolor=blue,citecolor=blue]{hyperref}
\usepackage{cite}
\usepackage[ruled]{algorithm2e}
\usepackage{fancybox}
\usepackage{framed}
\usepackage{enumitem}

\usepackage{amsthm}
\usepackage{mathtools}

\newtheorem{remark}{Remark}

\newtheorem{lemma}{\textbf{Lemma}}

\renewcommand\footnoterule{\kern-3pt \hrule width 2in \kern 2.6pt}

\journal{Nowhere}

\begin{document}

\begin{frontmatter}



\title{\LARGE A gathering  of Barbalat's lemmas and their (unsung) cousins}


\author{Zhiyong Sun}
\address{Control Systems Group, Department of EE, TU Eindhoven, the Netherlands}


\begin{abstract}
This note  presents a summary and review of various conditions for signal convergences, based on Barbalat-like lemmas and their variations. 
\end{abstract}


\end{frontmatter}


\section{Definitions and notations}
Consider a function $f(t) : \mathbb{R}_{\geq 0} \to \mathbb R$ that is locally integrable. Given a fixed $p\in (0, \infty)$, we say that $f(t)$ belongs to the $L^p$ space if $\int_0^{\infty} |f(s)|^p ds < \infty$. We define the function $p$-norm as $\Vert f(t) \Vert_{L^p} = ( \int_0^{\infty} |f(s)|^p ds)^\frac{1}{p}$, with a fixed $p\in [1, \infty)$. The $L^\infty$ norm for a function $f(t)$ is defined as $\|f(t)\|_\infty = \text{sup}_t |f(t)|$. We say a function $f(t) \in L^\infty$ if $\text{sup}_t |f(t)| < \infty$. 

For a vector-valued function $f(t): \mathbb{R}_{\geq 0} \to \mathbb R^n$, we say that $f(t)$ belongs to $L^p_n$ space (for a fixed $p\in (0, \infty)$) if $\int_0^{\infty} \Vert f(s)\Vert ^p  ds < \infty$ where $\Vert \cdot \Vert$ is the Euclidean norm\footnote{Note that any other vector norm could be used for the definition, due to equivalence of vector norms.}. We define the function $p$-norm in $\mathbb R^n$ as $\Vert f(t) \Vert_{L^p_n} = ( \int_0^{\infty} \Vert f(s)\Vert^p ds)^\frac{1}{p}$, with a fixed $p\in [1, \infty)$.
    
\section{Barbalat's lemma for signal convergence} 
The original version of the Barbalat's lemma is stated as below. 
\begin{lemma} \citep{khalil1996nonlinear}
Consider a scalar uniformly continuous function $f(t)$   such
that $\lim_{t \rightarrow \infty} \int_0^t f(\tau) d \tau$ exists and is finite. Then it holds $\lim_{t \rightarrow \infty}  f(t) = 0$.
\end{lemma} 

A sufficient condition to ensure that a function $f(t)$ is   uniformly continuous is that $\dot f(t)$ is bounded. Therefore, the first variation of the  Barbalat's lemma is stated as below. 
\begin{lemma}  
Consider a scalar  function $f(t)$   such
that $f(t) \in L^1$ and $\dot f(t) \in L^\infty$. Then it holds $\lim_{t \rightarrow \infty}  f(t) = 0$.
\end{lemma} 

\section{Cousins of Barbalat's lemma} 
Tao \citep{tao1997simple} proved the following simple alternative of Barbalat's lemma.

\begin{lemma} \citep{tao1997simple}
If $f(t) \in L^2$ and $\dot f(t) \in L^\infty$, then $\lim_{t \rightarrow \infty}  f(t) = 0$. 
\end{lemma}
\begin{remark}Some remarks are in order.
\begin{itemize}
 \item The condition ``$\dot f(t) \in L^\infty$" in the above lemma can be replaced by the more general condition that ``$f(t)$ is uniformly continuous"; indeed, the condition ``$\dot f(t) \in L^\infty$" is a sufficient condition to ensure that ``$f(t)$ is uniformly continuous".
 \item The proof of the above lemma also indicates that $f(t) \in L^\infty$, and therefore $f(t) \in L^p \cap L^\infty$, for any $p \in [2, \infty)$.

\item Rather than stating the condition ``$f(t) \in L^2$", a common condition (in many adaptive control books) is that ``$f(t) \in L^2 \cap L^\infty$", but the boundedness condition $f(t) \in L^\infty$ is not explicitly stated as it can be inferred by the two conditions on $\dot f(t)$ and $f(t)$ in the statement. 
\end{itemize}
\end{remark}

The following lemma is shown in \citep{desoer2009feedback}. 
\begin{lemma}\citep{desoer2009feedback}
Consider a scalar  function $f(t)$   such
that $f(t) \in L^2$ and $\dot f(t) \in L^2$. Then it holds that $f(t) \in  L^\infty$ and $\lim_{t \rightarrow \infty}  f(t) = 0$.
\end{lemma}

The following variation of Barbalat's lemma is given in \citep{krstic1995nonlinear}. 
\begin{lemma}  \citep[Page 491, Corollary A.7]{krstic1995nonlinear} 
    Consider the function   $f(t) : \mathbb{R}_{\geq 0} \to \mathbb R$. If $f, \dot f \in L^\infty$, and $f \in L^p$ for some $p \in [1. \infty)$, then $\lim_{t \rightarrow \infty}  f(t) = 0$. 
\end{lemma}

An alternative (and more general) lemma is proved in Tao's book \citep{tao2003adaptive}.
\begin{lemma} \citep[Page 80, Lemma 2.15]{tao2003adaptive}  \label{lemma:tao_barbalat}
 If $f(t) \in L^p$, $0 <p <\infty$ and $\dot f(t) \in L^\infty$, then $\lim_{t \rightarrow \infty}  f(t) = 0$. 
\end{lemma}
\begin{remark}
Note that in the statement of Lemma~\ref{lemma:tao_barbalat}, the condition  on the function $L^p$ space is relaxed to include $p 
\in (0, \infty)$. This is in contrast with the conditions in other Barbalat-like lemmas, while often one imposes the condition $p 
\in [1, \infty)$. 
\end{remark}
\section{A more general version of Barbalat-like lemma }
 
The following more general version of a Barbalat-like lemma is proved in \citep{farkas2016variations}.
\begin{lemma} \label{general:farkas}
 If $f(t) \in L^p$, $p \in [1, \infty)$ and $\dot f(t) \in L^q$, $q \in (1, \infty]$, then $f(t)$ is bounded (i.e., $f(t) \in L^\infty$) and uniformly continuous. Furthermore,   $\lim_{t \rightarrow \infty}  f(t) = 0$. 
\end{lemma}
\begin{remark}
Some remarks are in order.
\begin{itemize}
    \item In the above lemmas,  the conditions ``$\dot f(t) \in L^\infty$" or ``$\dot f(t) \in L^q$"  can be replaced by that ``$f(t)$ is uniformly continuous";
    \item The above lemmas for scalar function $f(t)$ can be extended to vector-valued functions $f(t) : \mathbb{R}_{\geq 0} \to \mathbb R^n$.  
    \item In Lemma~\ref{general:farkas}, the condition $\dot f(t) \in L^q$, $q \in (1, \infty]$ can be generalized to that $q \in [1, \infty]$ (i.e., the condition on $q$ also includes the case that $q=1$). Note: one can prove that if $\dot f(t) \in   L^1$ then $\lim_{t \rightarrow \infty}  f(t) $ exists and is finite. This fact together with the condition $f(t) \in L^p$, $p \in [1, \infty)$ will imply that $\lim_{t \rightarrow \infty}  f(t) = 0$. 
    \item The paper \citep{farkas2016variations} also discussed the rate of convergence of $L^p$ functions. Examples are given in \citep{farkas2016variations} to show the relations of different versions of Barbalat-like lemmas. 
    \item In the statement of Lemma~\ref{lemma:tao_barbalat}, the condition  on the function $L^p$ space is relaxed to include $p 
\in (0, \infty)$. It is possible to further relax the conditions of Lemma~\ref{general:farkas} that incorporate the conditions of Lemma~\ref{lemma:tao_barbalat}. 
\end{itemize}
\end{remark}

\section{Other versions of Barbalat-like lemmas with function compositions}
\begin{lemma} (See e.g., \citep{teel1999asymptotic})  \label{lemma:teel}
Let $\alpha(t): \mathbb{R}_{\geq 0} \to \mathbb{R}_{\geq 0}$ be   continuous, zero only at zero, and nondecreasing. If $f(t)$ is uniformly continuous on $[0, \infty)$ and $\alpha (f(t)) \in L^1$, then $\lim_{t \rightarrow \infty}  f(t) = 0$. 
\end{lemma}
 
A similar lemma is discussed in \citep{hou2010new}. 
\begin{lemma}
Let $M(z)$ be a continuous positive definite function  defined on $\{z: z \in \mathbb{R}^n, \|z\| \leq r\}$ for some $r$. If   $f(t): \mathbb{R}_{\geq 0} \to \mathbb{R}^n$ is a uniformly continuous function such that $f(t)\in L^\infty_n$ and $M(f(t)) \in L^1$, then $\lim_{t \rightarrow \infty}  f(t) = 0$.  
\end{lemma}
\begin{remark}
Some remarks are in order. 
\begin{itemize} 
\item The convergence result for scalar functions $f(t)$ in Lemma~\ref{lemma:teel} can be extended to vector-valued functions $f(t) : \mathbb{R}_{\geq 0} \to \mathbb R^n$. 
    \item By invoking the conditions of $L^p$ functions in Lemma~\ref{lemma:tao_barbalat} and Lemma~\ref{general:farkas}, it is possible to further relax the conditions of $L^1$ (i.e., $p=1$) function to other $L^p$ functions. 
    \item The two lemmas are often applied  (together with Lyapunov argument and comparison functions) to derive asymptotic convergence and $L^p$ stability of dynamical control systems. 
\end{itemize}
\end{remark}

\newpage
\bibliography{barbalat}
\bibliographystyle{elsarticle-harv}

\end{document}